\documentclass{article}

\usepackage{amsmath,amssymb,amsthm}
\usepackage{amsrefs}

\title{Weak Conservation Laws for Minimizers\\
which are not Pontryagin Extremals\footnote{Research report CM05/I-11.
Accepted for presentation (Paper No: 113) at the 2nd International Conference
``Physics and Control'' (PhysCon 2005), August 24-26, 2005, Saint Petersburg,
Russia. To appear in the respective Conference Proceedings.}}

\author{Delfim F.~M.~Torres\\
        \texttt{delfim@mat.ua.pt}}

\date{Department of Mathematics\\
      University of Aveiro\\
      3810-193 Aveiro, Portugal}

\newtheorem{theorem}{Theorem}

\theoremstyle{remark}

\theoremstyle{definition}
\newtheorem{definition}{Definition}

\begin{document}

\maketitle

\begin{abstract}
We prove a Noether-type symmetry theorem
for invariant optimal control problems with unrestricted controls.
The result establishes weak conservation laws along all the minimizers
of the problems, including those minimizers which do not satisfy
the Pontryagin Maximum Principle.
\end{abstract}

\smallskip


\smallskip


\smallskip



\section{Introduction}

Emmy Noether's first theorem \cite{MR0406752} is one of the most beautiful
and fundamental results of the calculus of variations.
The result comprises a universal principle, connecting the existence
of a family of transformations under which the functional
to be minimized is invariant (the existence of variational
symmetries) with the existence of conservation laws
(first integrals of the Euler-Lagrange differential equations).
Conservation laws can then be used to simplify the problem
of finding the minimizers. They have played an important role,
both in mathematics and physics, since the birth of the calculus
of variations in the eighteen century, having been extensively used
by giants like the Bernoulli brothers, Newton, Leibniz, Euler, Lagrange,
and Legendre. Conservation laws,
obtained from Noether's theorem, have a profound effect
on a vast number of disciplines, ranging from classical mechanics,
where they find important interpretations
such as conservation of energy, conservation of momentum,
or conservation of angular momentum, to engineering, economics,
control theory and their applications \cite{Gugushvili-et-al}.

The first extension of Noether's theorem to the more general
context of optimal control was published in 1973 \cite{MR0341229}.
Since then, many Noether-like theorems have been obtained in the context
of optimal control -- see \cite{MR2040245} and references therein.
We recall that all such versions of Noether's theorem
assume the Pontryagin maximum principle \cite{MR898009} to be satisfied,
and use its conditions, including the adjoint system, in their proofs.

Optimal control with unbounded controls is an area of strong current activity,
because of numerous applications involving modern technology such
as ``smart materials'' \cite{tr:Lasiecka:04}.
When there are no restrictions on the values of the control variables,
as in the calculus of variations,
it is well known that optimal control problems
may present solutions for which the Pontryagin Maximum Principle fails
to be satisfied (see \textrm{e.g.} \cite[\S 11.1]{MR1756410}).
This is due to the fact that the hypotheses of the existence theory
need to be complemented with additional regularity conditions
in order to proceed with the arguments which lead
to the maximum principle \cite{MR2006825}:
unboundedness of the controls ``propagates'' through the
dynamical control system, often causing a lack
of regularity for the solutions.

In spite of the gap between the hypotheses of
necessary optimality conditions and existence theorems,
J.~Ball proved \cite{MR755716} that, for time-invariant problems,
the conservation of the Hamiltonian (conservation of energy)
is still valid for minimizers which might
not satisfy the Euler-Lagrange necessary condition.
More recently, in 2002, G.~Francfort and J.~Sivaloganathan
proposed a generalization of Ball's result,
giving some applications to hyper-elasticity \cite{MR1950811}.

In this note we extend the previous results \cite{MR1950811}
from the calculus of variations framework to the
optimal control setting. We obtain weak conservation
laws for minimizers which do not necessarily satisfy
the Pontryagin maximum principle.


\section{Bad Behavior in Optimal Control}

The optimal control problem consists to minimize a cost functional
\begin{equation*}
I\left[x(\cdot),u(\cdot)\right] = \int_a^b L\left(t,x(t),u(t)\right) dt
\end{equation*}
subject to a control system described by ordinary differential equations
\begin{equation}
\label{eq:tr:prP:cs}
\dot{x}(t) = \varphi\left(t,x(t),u(t)\right)
\end{equation}
together with certain appropriate endpoint conditions.
The Lagrangian $L : [a,b] \times \mathbb{R}^n \times \mathbb{R}^m
\rightarrow \mathbb{R}$ and the velocity vector $\varphi : [a,b]
\times \mathbb{R}^n \times \mathbb{R}^m \rightarrow \mathbb{R}^n$
are given, and assumed to be smooth: $L(\cdot,\cdot,\cdot)$,
$\varphi(\cdot,\cdot,\cdot)$ $\in C^1$.
We are interested in the case where there are no restrictions
on the control set: $u(t) \in \mathbb{R}^m$. We denote the problem by $(P)$.
In the particular case $\varphi(t,x,u) = u$, one obtains the fundamental problem
of the calculus of variations, which covers all classical mechanics.
The choice of the classes $\mathcal{X}$ and $\mathcal{U}$, respectively
of the state $x : [a,b] \rightarrow \mathbb{R}^n$ and control variables
$u : [a,b] \rightarrow \mathbb{R}^m$, play an important role in our discussion.
In connection with the Pontryagin maximum principle, the optimal controls
are typically assumed to be essentially bounded, $\mathcal{U} = L_{\infty}$
\cite{MR898009}; while to guarantee existence, compactness arguments require
a bigger class of measurable control functions, $\mathcal{U} = L_{1}$
\cite{MR688142}. Given an optimal control problem with unrestricted controls,
it may happen that the Pontryagin maximum principle is valid,
while existence of minimizers is not guaranteed;
or it may happen that the minimizers predicted by the
existence theory fail to be Pontryagin extremals.
Conservation laws are obtained from Noether's theorem,
assuming that the minimizers are Pontryagin extremals.
In this work our objective is to prove weak conservation laws
valid for bad-behaved problems with minimizers which
are not Pontryagin extremals. We begin to explain
why bad-behavior can occur.

The Pontryagin maximum principle is a necessary optimality
condition which can be obtained from a general Lagrange multiplier
theorem in spaces of infinite dimension
(\textrm{cf.} \textrm{e.g.} \cite{MR1846384,MR2014219}).
Introducing the Hamiltonian function
\begin{equation}
\label{eq:tr:H}
H(t,x,u,\psi) = - L(t,x,u) + \psi \cdot \varphi(t,x,u) \, ,
\end{equation}
where $\psi_i$, $i = 1,\ldots,n$, are the ``Lagrange multipliers''
or the ``generalized momenta'', the multiplier theorem asserts that
the optimal control problem is equivalent to the maximization
of the augmented functional
\begin{equation}
\label{eq:tr:af}
J\left[x(\cdot),u(\cdot),\psi(\cdot)\right] = \int_a^b \left(
H\left(t,x(t),u(t),\psi(t)\right) - \psi(t) \cdot \dot{x}(t) \right) dt \, .
\end{equation}
Let us assume, for simplicity,
$\mathcal{X} = C^1\left([a,b];\mathbb{R}^n\right)$,
$\mathcal{U} = C\left([a,b];\mathbb{R}^m\right)$.
Let $\left(\tilde{x}(\cdot),\tilde{u}(\cdot),\tilde{\psi}(\cdot)\right)$
solve the problem, and consider arbitrary $C^1$-functions
$h_1 \, , h_3 : [a,b] \rightarrow \mathbb{R}^n$,
$h_1(\cdot)$ vanishing at $a$ and $b$ ($h_1(\cdot) \in C^1_0\left([a,b]\right)$),
and arbitrary continuous $h_2 : [a,b] \rightarrow \mathbb{R}^m$.
Let $\varepsilon$ be a scalar. By definition of maximizer, we have
\begin{equation*}
J\left[(\tilde{x}+\varepsilon h_1)(\cdot),(\tilde{u}+\varepsilon h_2)(\cdot),
(\tilde{\psi}+\varepsilon h_3)(\cdot)\right] \le
J\left[\tilde{x}(\cdot),\tilde{u}(\cdot),\tilde{\psi}(\cdot)\right] \, ,
\end{equation*}
and one has the following necessary condition:
\begin{equation}
\label{eq:tr:ncci}
\frac{d}{d\varepsilon} \left.
J\left[(\tilde{x}+\varepsilon h_1)(\cdot),(\tilde{u}+\varepsilon h_2)(\cdot),(\tilde{\psi}
+\varepsilon h_3)(\cdot)\right]
\right|_{\varepsilon = 0} = 0\, .
\end{equation}
Differentiating \eqref{eq:tr:ncci} gives
\begin{multline}
\label{eq:tr:nc:bip}
0 = \int_a^b \left[ \frac{\partial H}{\partial x}\left(t,\tilde{x}(t),
\tilde{u}(t),\tilde{\psi}(t)\right) \cdot h_1(t)
+ \frac{\partial H}{\partial u}\left(t,\tilde{x}(t),
\tilde{u}(t),\tilde{\psi}(t)\right) \cdot h_2(t) \right. \\
\left. + \frac{\partial H}{\partial \psi}\left(t,\tilde{x}(t),
\tilde{u}(t),\tilde{\psi}(t)\right) \cdot h_3(t)
- h_3(t) \cdot \dot{\tilde{x}}(t) - \tilde{\psi}(t) \cdot \dot{h}_1(t)\right] dt \, .
\end{multline}
Integrating the $\tilde{\psi}(t) \cdot \dot{h}_1(t)$ term by parts,
and having in mind that $h_1(a) = h_1(b) = 0$, one derives
\begin{multline}
\label{eq:tr:afterPP}
\int_a^b \left[ \left(\frac{\partial H}{\partial x}\left(t,\tilde{x}(t),
\tilde{u}(t),\tilde{\psi}(t)\right) + \dot{\tilde{\psi}}(t)\right) \cdot h_1(t) \right.
+ \frac{\partial H}{\partial u}\left(t,\tilde{x}(t),
\tilde{u}(t),\tilde{\psi}(t)\right) \cdot h_2(t) \\
\left. + \left(\frac{\partial H}{\partial \psi}\left(t,\tilde{x}(t),
\tilde{u}(t),\tilde{\psi}(t)\right) - \dot{\tilde{x}}(t)\right) \cdot h_3(t)\right] dt = 0 \, .
\end{multline}
Note that \eqref{eq:tr:afterPP} was obtained for any variation $h_1(\cdot)$,
$h_2(\cdot)$, and $h_3(\cdot)$. Choosing $h_1(t) = h_2(t) \equiv 0$, and
$h_3(\cdot)$ arbitrary, one obtains the control system
\eqref{eq:tr:prP:cs}:
\begin{equation}
\label{eq:tr:PMP:cs}
\dot{\tilde{x}}(t) = \frac{\partial H}{\partial \psi}\left(t,\tilde{x}(t),
\tilde{u}(t),\tilde{\psi}(t)\right) \, , \quad t \in [a,b] \, .
\end{equation}
With $h_1(\cdot)$ arbitrary, and $h_2(t) = h_3(t) \equiv 0$,
we obtain the \emph{adjoint system}:
\begin{equation}
\label{eq:tr:PMP:as}
\dot{\tilde{\psi}}(t) = - \frac{\partial H}{\partial x}\left(t,\tilde{x}(t),
\tilde{u}(t),\tilde{\psi}(t)\right) \, , \quad t \in [a,b] \, .
\end{equation}
Finally, with $h_2(\cdot)$ arbitrary, and $h_1(t) = h_3(t) \equiv 0$,
the \emph{maximality condition} is obtained:
\begin{equation}
\label{eq:tr:PMP:mc}
\frac{\partial H}{\partial u}\left(t,\tilde{x}(t),
\tilde{u}(t),\tilde{\psi}(t)\right) = 0 \, , \quad t \in [a,b] \, .
\end{equation}
A necessary optimality condition for $\left(\tilde{x}(\cdot),\tilde{u}(\cdot)\right)$
to be a minimizer of problem $(P)$ is given by the \emph{Pontryagin maximum principle}:
there exists $\tilde{\psi}(\cdot)$ such that the 3-tuple
$\left(\tilde{x}(\cdot),\tilde{u}(\cdot),\tilde{\psi}(\cdot)\right)$ satisfy
all the conditions \eqref{eq:tr:PMP:cs}, \eqref{eq:tr:PMP:as}, and \eqref{eq:tr:PMP:mc}.
We recall that conditions \eqref{eq:tr:PMP:cs}, \eqref{eq:tr:PMP:as},
and \eqref{eq:tr:PMP:mc} imply the equality
\begin{equation}
\label{eq:tr:PMP:dHdt}
\frac{d}{dt} H\left(t,\tilde{x}(t),\tilde{u}(t),\tilde{\psi}(t)\right)
= \frac{\partial H}{\partial t}\left(t,\tilde{x}(t),\tilde{u}(t),\tilde{\psi}(t)\right) \, .
\end{equation}
For piecewise smooth state trajectories, $\mathcal{X} = PC^1\left([a,b];\mathbb{R}^n\right)$,
and piecewise continuous controls, $\mathcal{U} = PC\left([a,b];\mathbb{R}^m\right)$,
similar arguments than those used to derive conditions \eqref{eq:tr:PMP:cs}, \eqref{eq:tr:PMP:as},
and \eqref{eq:tr:PMP:mc} continue to be justifiable. In fact, as already mentioned,
the arguments can be carried out for essentially bounded controls,
$\mathcal{U} = L_{\infty}\left([a,b];\mathbb{R}^m\right)$,
and Lipschitzian state variables, $\mathcal{X} = W_{1,\infty}\left([a,b];\mathbb{R}^n\right)$.
But if one took $\mathcal{U} = L_1\left([a,b];\mathbb{R}^m\right)$,
and $\mathcal{X} = W_{1,1}\left([a,b];\mathbb{R}^n\right)$,
as required by the existence theory, integration by parts of
$\tilde{\psi}(t) \cdot \dot{h}_1(t)$ in \eqref{eq:tr:nc:bip}
can no longer be justified, and one can not conclude with the adjoint system \eqref{eq:tr:PMP:as}
and equality \eqref{eq:tr:PMP:dHdt}.
This is more than a technical difficulty, and explains the possibility of bad-behavior
illustrated by the Ball-Mizel example \cite{MR801585}.
In the Ball-Mizel problem one has $n=m=1$,
$L(t,x,u)=(x^3 - t^2)^2 u^{14} + \epsilon u^2$,
and $\varphi(t,x,u)=u$:
\begin{gather*}
\int_{0}^{1}\left[\left(x^{3}(t)-t^{2}\right)^{2}\,u(t)^{14}
+\epsilon \,u(t)^{2}\right]\,dt\longrightarrow \min \, ,\\
\dot{x}\left(t\right)=u\left(t\right) \, ,\\
x\left(0\right)=0 \, , \quad x\left(1\right) = k \, .
\end{gather*}
For some values of the constants $\epsilon$ and $k$
there exists the unique optimal control
$\tilde{u}(t) =\frac{2k}{3}\,t^{-1/3}$ \cite{MR756513},
which belongs to $L_1$ but not to $L_\infty$. The Pontryagin maximum principle
is not satisfied since the adjoint system \eqref{eq:tr:PMP:as}
\begin{equation*}
\dot{\tilde{\psi}}(t)=-\frac{\partial H}{\partial x}\left(t,\tilde{x}(t),
\tilde{u}(t),\tilde{\psi}(t)\right) \cong t^{-4/3}
\end{equation*}
is not integrable. In this paper we obtain a new version
of Noether's theorem, without using the adjoint system
\eqref{eq:tr:PMP:as} and the property \eqref{eq:tr:PMP:dHdt}.
This makes Noether's principle valid both for well- and bad-behaved
optimal control problems.


\section{Conservation Laws in Optimal Control}

In 1918 Emmy Noether established the key result to find conservation
laws in the calculus of variations \cite{MR0406752,MR2004181}. We sketch here
the standard argument used to derive Noether's theorem and conservation laws
in the optimal control setting (\textrm{cf.} \textrm{e.g.} \cite{MR0341229,MR1901565}).

Let us consider a one-parameter group of $C^1$-transformations of the form
\begin{equation}
\label{eq:tr:pt}
h^s(t,x,u,\psi) = \left(h_t(t,x,u,\psi,s),
h_x(t,x,u,\psi,s),h_u(t,x,u,\psi,s),h_\psi(t,x,u,\psi,s)\right) \, ,
\end{equation}
where $s$ denote the independent parameter of the transformations.
We require that to the parameter value $s = 0$ there corresponds
the identity transformation:
\begin{equation}
\label{eq:tr:pts0}
\begin{split}
h^0(t,x,u,\psi) &= \left(h_t(t,x,u,\psi,0),
h_x(t,x,u,\psi,0),h_u(t,x,u,\psi,0),h_\psi(t,x,u,\psi,0)\right) \\
&= (t,x,u,\psi) \, .
\end{split}
\end{equation}
Associated to the group of transformations \eqref{eq:tr:pt} we consider
the \emph{infinitesimal generators}
\begin{gather}
T(t,x,u,\psi) = \left.\frac{d}{ds} h_t(t,x,u,\psi,s)\right|_{s = 0} \, , \quad
X(t,x,u,\psi) = \left.\frac{d}{ds} h_x(t,x,u,\psi,s)\right|_{s = 0} \, , \label{eq:tr:TX} \\
U(t,x,u,\psi) = \left.\frac{d}{ds} h_u(t,x,u,\psi,s)\right|_{s = 0} \, , \quad
\Psi(t,x,u,\psi) = \left.\frac{d}{ds} h_\psi(t,x,u,\psi,s)\right|_{s = 0} \, . \notag
\end{gather}

\begin{definition}
\label{def:tr:inv}
The optimal control problem $(P)$ is said to be invariant under
a one-parameter group of $C^1$-transformations \eqref{eq:tr:pt}
if, and only if,
\begin{multline}
\label{eq:tr:inv}
\left.\frac{d}{ds}
  \Biggl\{
    \biggl[
      H\left(h^s\left(t,x(t),u(t),\psi(t)\right)\right) \right. \\
      \left. - h_\psi\left(t,x(t),u(t),\psi(t),s\right) \cdot
        \frac{\frac{dh_x\left(t,x(t),u(t),\psi(t),s\right)}{dt}}{\frac{dh_t\left(t,x(t),u(t),\psi(t),s\right)}{dt}}
    \biggr]
    \frac{dh_t\left(t,x(t),u(t),\psi(t),s\right)}{dt}
  \Biggr\}
\right|_{s = 0} = 0 \, ,
\end{multline}
with $H$ the Hamiltonian \eqref{eq:tr:H}.
\end{definition}

Having in mind \eqref{eq:tr:pts0}, condition \eqref{eq:tr:inv}
is equivalent to
\begin{multline}
\label{eq:tr:cond:nec:suf:inv}
\frac{\partial H}{\partial t} T +
\frac{\partial H}{\partial x} \cdot X + \frac{\partial H}{\partial u} \cdot U +
\frac{\partial H}{\partial \psi} \cdot \Psi - \Psi \cdot \dot{x}(t)
- \psi(t) \cdot \frac{d}{dt}X + H  \frac{d}{dt} T = 0 \, ,
\end{multline}
where here, and to the end of the paper, all functions are evaluated
at $\left(t,x(t),u(t),\psi(t)\right)$ whenever not indicated.
Along a Pontryagin extremal $\left(x(\cdot),u(\cdot),\psi(\cdot)\right)$
equalities \eqref{eq:tr:PMP:cs}, \eqref{eq:tr:PMP:as}, \eqref{eq:tr:PMP:mc},
and \eqref{eq:tr:PMP:dHdt} are in force, and \eqref{eq:tr:cond:nec:suf:inv}
reduces to
\begin{equation*}
\frac{dH}{dt} T - \dot{\psi}(t) \cdot X
- \psi(t) \cdot \frac{dX}{dt} + H \frac{dT}{dt} = 0
\Leftrightarrow \frac{d}{dt} \left(
\psi(t) \cdot X - H T\right) = 0 \, .
\end{equation*}
We have just proved \emph{Noether's theorem for optimal control problems}.
\begin{theorem}[Noether's Theorem]
If the optimal control problem is invariant under \eqref{eq:tr:pt},
in the sense of Definition~\ref{def:tr:inv}, then
\begin{multline}
\label{eq:tr:cl}
\psi(t) \cdot X\left(t,x(t),u(t),\psi(t)\right)
- H\left(t,x(t),u(t),\psi(t)\right) T\left(t,x(t),u(t),\psi(t)\right) = c
\end{multline}
($c$ a constant; $t \in [a,b]$; $T$ and $X$ are given according to \eqref{eq:tr:TX};
$H$ is the Hamiltonian \eqref{eq:tr:H}) is a \emph{conservation law},
that is, \eqref{eq:tr:cl} is valid along all the minimizers
$(x(\cdot),u(\cdot))$ of $(P)$ which are Pontryagin extremals.
\end{theorem}

All available versions of Noether's theorem found in the literature
are valid only for well-behaved optimal control problems
(conservation laws are, by definition, valid for minimizers which are
Pontryagin extremals). In the next section we provide
the first optimal control version of Noether's theorem valid
in presence of bad behavior (valid also for minimizers which
are not Pontryagin extremals). For that we need a new notion
of conservation law.


\section{Weak Conservation Laws in Optimal Control}

In 1879 Paul duBois-Reymond proved an important basic result.
From duBois-Reymond lemma we know that
\begin{equation*}
\int_a^b \left(
\psi(t) \cdot X - H T\right) \dot{\theta}(t) dt = 0 \, ,
\quad \forall \, \theta(\cdot) \in C_0^1([a,b]),
\end{equation*}
is a weak form of conservation law \eqref{eq:tr:cl}.
Follows the main result of the paper.
\begin{theorem}
\label{MainResult:wcl:torres}
If the optimal control problem is invariant under \eqref{eq:tr:pt},
in the sense of Definition~\ref{def:tr:inv}, then
\begin{multline*}
\int_a^b \biggl[\psi(t) \cdot X\left(t,x(t),u(t),\psi(t)\right) \\
- H\left(t,x(t),u(t),\psi(t)\right) T\left(t,x(t),u(t),\psi(t)\right)\biggr]
\dot{\theta}(t) dt = 0
\end{multline*}
($\theta(\cdot)$ is an arbitrary $W_{1,1}([a,b]; \mathbb{R})$ function satisfying
$\theta(a) = \theta(b) = 0$) holds along all the minimizers
$(x(\cdot),u(\cdot),\psi(\cdot)) \in W_{1,1} \times L_1 \times W_{1,1}$ of \eqref{eq:tr:af}.
\end{theorem}
\begin{proof}
Replacing the parameter $s$ of the group \eqref{eq:tr:pt} by function $s \theta(t)$,
the infinitesimal generators are then given by
\begin{equation*}
\begin{split}
\left.\frac{d}{ds} h_t\left(t,x,u,\psi,s\theta(t)\right)\right|_{s = 0}
    &= T(t,x,u,\psi) \theta(t)\, , \\
\left.\frac{d}{ds} h_x\left(t,x,u,\psi,s\theta(t)\right)\right|_{s = 0}
    &= X(t,x,u,\psi) \theta(t) \, ,  \\
\left.\frac{d}{ds} h_u\left(t,x,u,\psi,s\theta(t)\right)\right|_{s = 0}
    &= U(t,x,u,\psi) \theta(t)\, , \\
\left.\frac{d}{ds} h_\psi\left(t,x,u,\psi,s\theta(t)\right)\right|_{s = 0}
    &= \Psi(t,x,u,\psi) \theta(t) \, ,
\end{split}
\end{equation*}
with $T$, $X$, $U$, and $\Psi$ as in \eqref{eq:tr:TX},
and the necessary and sufficient condition of invariance
\eqref{eq:tr:cond:nec:suf:inv} takes the form
\begin{multline}
\label{eq:tr:cond:nec:suf:inv:2NT}
\left(\frac{\partial H}{\partial t} T +
\frac{\partial H}{\partial x} \cdot X + \frac{\partial H}{\partial u} \cdot U +
\frac{\partial H}{\partial \psi} \cdot \Psi - \Psi \cdot \dot{x}(t)
+ H  \frac{dT}{dt} - \psi(t) \cdot \frac{dX}{dt}\right) \theta(t) \\
+ \left(H T - \psi(t) \cdot X\right) \dot{\theta}(t) = 0 \, .
\end{multline}
Condition \eqref{eq:tr:nc:bip} with $h_1(t) = X \theta(t)$,
$h_2(t) = U \theta(t)$, and $h_3(t) = \Psi \theta(t)$, gives
\begin{multline}
\label{eq:tr:cond:nec:min:2NT}
\int_a^b \biggl[\left(
\frac{\partial H}{\partial x} \cdot X + \frac{\partial H}{\partial u} \cdot U +
\frac{\partial H}{\partial \psi} \cdot \Psi - \Psi \cdot \dot{x}(t)
- \psi(t) \cdot \frac{dX}{dt}\right) \theta(t)\\
-\left(\psi(t) \cdot X\right) \dot{\theta}(t)\biggr] dt = 0 \, .
\end{multline}
Using \eqref{eq:tr:cond:nec:suf:inv} in \eqref{eq:tr:cond:nec:min:2NT}
permits to write
\begin{equation}
\label{eq:tr:invmin2nta}
\int_a^b \biggl[\left(
-\frac{\partial H}{\partial t} T - H \frac{dT}{dt}\right) \theta(t)
-\left(\psi(t) \cdot X\right) \dot{\theta}(t)\biggr] dt = 0 \, ;
\end{equation}
while, on the other hand, using \eqref{eq:tr:cond:nec:suf:inv:2NT}
in \eqref{eq:tr:cond:nec:min:2NT}, one obtains
\begin{equation}
\label{eq:tr:invmin2ntb}
\int_a^b \biggl[\left(
\frac{\partial H}{\partial t} T + H \frac{dT}{dt}\right) \theta(t)
+\left(H T\right) \dot{\theta}(t)\biggr] dt = 0 \, .
\end{equation}
The conclusion follows summing up \eqref{eq:tr:invmin2nta} and \eqref{eq:tr:invmin2ntb}:
\begin{equation*}
\int_a^b \left(H T - \psi(t) \cdot X\right) \dot{\theta}(t) dt = 0 \, .
\end{equation*}
\end{proof}

Invariance under an infinite continuous group of transformations, which
rather than dependence on parameters depend upon arbitrary functions,
is considered by Noether in the original paper \cite{MR0406752}.
This is sometimes called ``the second Noether theorem''.
We refer the reader to \cite{MR1980565} for an extension of the second Noether
theorem to optimal control problems which are semi-invariant under symmetries
depending upon $k$ arbitrary functions of the independent variable and their derivatives.
Theorem~\ref{MainResult:wcl:torres} is easily formulated under more general
notions of invariance.

Theorem~\ref{MainResult:wcl:torres} gives, from the invariance properties
of the optimal control problems, weak conservation laws along all the minimizers,
including those minimizers which does not satisfy the standard
necessary optimality conditions like the Pontryagin Maximum
Principle or the Euler-Lagrange differential equations. Such fact
may be useful to identify more general classes of well-
and bad-behaved problems in the calculus of variations and optimal
control, \textrm{e.g.}, to synthesize a broad class of invariant
problems exhibiting the Lavrentiev phenomenon \cite{MR2043868}.
This possibility is under investigation and will be addressed elsewhere.


\section*{Acknowledgements}

This work was partially supported by the
Portuguese Foundation for Science and Technology
(FCT) through the Control Theory Group (cotg)
of the Centre for Research in Optimization and Control (CEOC).



\end{document}